\documentclass[psfig]{article}

\usepackage{amssymb, amsmath, amsthm}

\newtheorem{theo}{Theorem}[section]
\newtheorem{lem}{Lemma}[section]

\newcommand{\ov}{\overline}

\newcommand{\BB}[1]{\ensuremath{\mathbb{#1}}}

\newcommand{\Q}{\ensuremath{\BB{Q}}}
\newcommand{\R}{\ensuremath{\BB{R}}}
\newcommand{\Z}{\ensuremath{\BB{Z}}}

\newcommand{\be}{\begin{equation}}
\newcommand{\ee}{\end{equation}}
\newcommand{\bes}{\begin{equation*}}
\newcommand{\ees}{\end{equation*}}
\newcommand{\bi}{\begin{itemize}}
\newcommand{\ei}{\end{itemize}}
\newcommand{\bea}{\begin{eqnarray}}
\newcommand{\eea}{\end{eqnarray}}
\newcommand{\beas}{\begin{eqnarray*}}
\newcommand{\eeas}{\end{eqnarray*}}

\begin{document}
\title{Multiplicity and regularity of large periodic solutions with rational frequency for a class of semilinear monotone wave equations.}
\author{Jean-Marcel Fokam \\
{\small \tt fokam\,@aun.edu.ng}\\
School of Arts and Sciences, American University of Nigeria\\
Yola, 
}


\date{} 



\maketitle

\begin{abstract}
Nous {$\rm{d\acute{e}montrons}$} l'existence d'une ${\rm{infinit\acute{e}}}$ de solutions fortes, de norme grande, pour une classe d'{$\rm{\acute{e}quations}$} {$\rm{semilin\acute{e}aires}$} avec des conditions {$\rm{p\acute{e}riodiques}$} sur le bord:\\
\[
u_{tt}-u_{xx}=f(x,u),
\]
\[
 u(0,t)=u(\pi,t)\,\,\,, u_x(0,t)=u_x(\pi,t).
\]
Notre {$\rm{m\acute{e}thode}$} est bas${\rm{\grave{e}}}$e sur de nouvelles estimations pour le {$\rm{probl\grave{e}me}$} {$\rm{lin\acute{e}aire}$} avec conditions {$\rm{p\acute{e}riodiques}$} sur le bord, en combinant les {$\rm{m\acute{e}thodes}$} de Littlewood-Paley, le {$\rm{th\acute{e}or\grave{e}me}$} de Hausdorff-Young et une formulation variationelle de Rabinowitz, \cite{R78},\cite{R84}. Nous contruisons une nouvelle approche pour la {$\rm{r\acute{e}gularite}$} des solutions au sens des distributions en {$\rm{d\acute{e}rivant}$} les {$\rm{\acute{e}quations}$} et en utilisant les estimations de type Gagliardo-Nirenberg. \\
We prove the existence of infinitely many classical large periodic solutions for a class of semilinear wave equations with periodic boundary conditions:
\[
u_{tt}-u_{xx}=f(x,u),
\]
\[
 u(0,t)=u(\pi,t)\,\,\,, u_x(0,t)=u_x(\pi,t).
\]
Our argument relies on some new estimates for the linear problem with periodic boundary conditions, by combining Littlewood-Paley techniques, the Hausdorff-Young theorem of harmonic analysis, and a variational formulation due to Rabinowitz \cite{R78},\cite{R84}. We also develop a new approach to the regularity of the distributional solutions by differentiating the equations and employing Gagliardo-Nirenberg estimates.\\

\end{abstract}

\newpage

\section{Introduction}

\noindent In this paper we construct infinitely many large classical time-periodic solutions
for the following semilinear wave equation: 
\begin{equation}
    u_{tt}-u_{xx}-f(x,u)=0 \label{papier}
\end{equation}
\be
    u(0,t)=u(\pi,t),\,\ u_x(0,t)=u_x(\pi,t)
\label{bdry}
\ee
where $f$ is $C^{2,1}$, has polynomial growth and depends on $x,u$. The existence of large periodic solutions with periodic boundary conditions is not well understood. As $u=0$ is a trivial solution we seek here nontrivial solutions of (\ref{papier}),(\ref{bdry}). When the frequency is irrational the method of Craig and Wayne in \cite{CW}, extended to higher dimension by Bourgain\cite{Bourgain95} and Berti and Bolle \cite{BertiBolle2010}proves the existence of small periodic solutions for typical potentials but the existence of classical periodic solutions for rational frequency is not known. Note that typical constant potentials in \cite{CW},\cite{Bourgain95}, \cite{BertiBolle2010} are satisfied for
\be
   u_{tt}-\Delta u-m-f(x,u)=0
\ee
for \textit{typical m} which exclude $m=0$. The so-called {\textit{resonant}} case $m=0$ and with $f(x,u)$ independent of $x$ for periodic boundary conditions has been studied by Berti and Procesi in \cite{BertiProcesi}. The lack of $x$ dependence in \cite{BertiProcesi} allows  to employ ordinary differential equations techniques and they showed the existence of quasi-periodic solutions where the frequency vector depends on two frequencies $(\omega_1,\omega_2(\epsilon))$. While they consider $\omega_1\in\Q$, their results do not imply the existence of periodic solutions with rational frequency as $\omega_2(\epsilon)$ there is never rational. Chierchia and You in \cite{ChierchiaYou} study the problem with periodic boundary conditions and a potential :
\be
    u_{tt}-u_{xx}-v(x)u-f(u)=0
\ee
where $f$ only depends on $u$, however their method excludes the constants potentials $v(x)=m$. Bricmont, Kupiainen and Schenkel in \cite{BKS} prove the existence of quasi-periodic solutions with periodic boundary conditions in the non-resonant case $m>0$ and $f$ depending only on $u$. In \cite{BKS} they find quasi-periodic solutions for a set of positive measure of frequencies hence prove the existence of quasi-periodic solutions for irrational frequencies.

On the other hand there exists a substantial amount of literature for semilinear wave equations with Dirichlet boundary conditions see for instance \cite{R84},\cite{R78}, \cite{BCN} for rational frequencies and the proofs of existence of \textit{classical solutions with $f$ having some spatial dependence rely on a fundamental solution discovered by Lovicarova} in \cite{Lovicarova}. The existence of periodic solutions, with irrational frequencies with Dirichlet boundary conditions in the resonant case ($m=0$) was shown by Lidskii and Schulman \cite{LidskiiShulman}, by Bambusi in \cite{Bambusi}, Bambusi and Paleari in \cite{BambusiPaleari}, Berti and Bolle \cite{BertiBolleDuke}, \cite{BertiBolle2003} and for quasi-periodic solutions in Yuan \cite{Yuan}. Quasi-periodic solutions with Dirichlet boundary conditions via KAM techniques has been shown by P\"oschel \cite{Poschel}, Kuksin \cite{Kuksin} and Wayne \cite{Wayne}. 

De Simon and Torelli in \cite{DST} do not employ Lovicarova's formula but their $C^0$ estimate relies on $L^2$ a priori estimates on $f(x,u)$ which are not readily avaliable for distributional solutions of (\ref{papier}). The difficulty in proving regularity of distributional solution of (\ref{papier}) stems for the kernel of $\Box$ which is infinite dimensional. In absence of a fundamental solution for the d'Alembertian under periodic boundary conditions problem we develop an approach based on tools from harmonic analysis such as Littlewood-Paley techniques, the Hausdorff-Young theorem and Gagliardo-Nirenberg estimates. The Hausdorff-Young theorem had been employed earlier by  Willem in \cite{Willem85} to get a $L^\infty$  a priori estimate on solutions, by Coron to prove a Sobolev embedding in \cite{Coron} and by Zhou in \cite{Zhou87}. The argument we give here to prove the Sobolev embedding in \cite{Coron} follows the Fourier approach to the Sobolev embedding as in the notes by Chemin, \cite{Chemin}. We do prove a stronger estimate than the one in \cite{Coron}, which provides information about the best constant of the Sobolev embedding. Our argument also shows that the Sobolev embedding in \cite{Coron} is \textit{compact}. In this paper the existence of \textit{classical solution for time periodic solutions with periodic boundary conditions of the semilinear wave equation (\ref{papier}) will be shown by proving the stronger $C^\gamma$ H\"older estimates than the $L^\infty$ in} \cite{Willem85}, and our approach also gives an alternative proof of the existence of classical periodic solutions in the case of Dirichlet boundary conditions with semilinear term with some spatial dependence for $f(x,u)$ sufficiently smooth in both arguments $x$ and $u$.
 
In section 1 we prove the linear estimates we need to prove the regularity of the solution. In section 2 we follow the scheme of \cite{R84} and \cite{R78} to construct weak solutions and in section 3 we show the regularity of the solution by repeated differentiation of the equations, the linear estimates proved in section 1 and Gagliardo-Nirenberg inequalities .

Since our proof is of variational nature it is natural to ask if there is a notion of critical exponent or critical growth for this equation. An open question is then whether there are semilinear terms $f(x,u)$ of say exponential or super exponential type (as this paper deals with semilinear terms of polynomial type) for which there are large amplitude distributional solutions which are not classical ($f(x,u)$ being assumed to be smooth).


We seek time-periodic solutions satisfying periodic boundary conditions so we seek functions $u\in\R$ with expansions of the form
\[
u(x,t)=\sum_{(j,k)\in\Z\times\Z}\widehat{u}(j,k)e^{i2jx}e^{ikt}
\]
and define the function space $E$:
\[
||u||^2_{E}=\sum_{2j\neq \pm k}\frac{|Q|}{4}|k^2-4j^2||\widehat{u}(j,k)|^2+\sum_{2j=\pm k}|4j^2||\widehat{u}(j,k)|^2+|\widehat{u}(0,0)|^2
\]
where $Q=[0,\pi]\times[0,2\pi]$ and define the functions spaces $E^+,E^-,N$ as follows:
\[
N=\{ u\in E, \,\, \widehat{u}(j,k)=0 \,\, {\rm{for}} \,\, 2|j|\neq |k|\}.
\]
Note that in the case of periodic boundary conditions the structure of the kernel $N$ of $\Box$ is slightly different than in the case of Dirichlet boundary conditions. Here $v\in N$ we have
\begin{eqnarray}v(x,t)
&  = & \sum_{j=\pm k}\widehat{v(j,k)}e^{i2jx+ikt}\nonumber\\
& =  & \sum_{j\neq 0}\widehat{v(j,2j)}e^{i2j(x+t)}+\sum_{j}\widehat{v(j,-2j)}e^{i2j(x-t)}
\end{eqnarray}
and
\be
v(x,t)=p_1(x+t)+p_2(x-t)=v^+(x,t)+v^-(x,t)
\ee
where $v^+(x,t)=p_1(x,t)$, $v^-(x,t)=p_2(x,t)$, where the $p_1,p_2\in H^1(0,\pi)$ $\pi$-periodic functions and defined as $p_1(s)=\sum_{j}\widehat{p_1(j)}e^{i2js}$, $p_2(s)=\sum_{j}\widehat{p_2(j)}e^{i2js}$ and $\widehat{p_1(0)}=0$, $\widehat{p_1(j)}=\widehat{v(j,2j)}$,$\widehat{p_2(j)}=\widehat{v(j,-2j)}$.
\[
E^+=\{ u\in E, \,\, \widehat{u}(j,k)=0 \,\, {\rm{for}} \,\, |k|\leq 2|j|\}
\]
\[
E^-=\{ u\in E, \,\, \widehat{u}(j,k)=0 \,\, {\rm{for}} \,\, |k|\geq 2|j|\},
\]
$u=v+w$, $w=w^++w^-$ where $w\in E$, $w^+\in E^+$,$w^-\in E^-$ and $v\in N$  
and define the norm on $E\oplus N$
\[
||u||_{\beta,E}^2=||w^+||^2_E+||w^-||^2_E+\beta||v||^2_{H^1}.
\]
\be
I_\beta(u)=\int_Q[\frac{1}{2}(u_t^2-u_x^2-\beta (v^2+v_t^2)-F(x,u)]dxdt.
\ee
When $u$ is trigonometric polynomial, $I_\beta$ can also be represented in $E^m\oplus N^m$ as:
\[
I_\beta(u)=\frac{1}{2}(||w^+||_E^2-||w^-||_E^2-\beta(||v||^2_{L^2}+||v_t||^2_{L^2})-\int_Q F(x,u)dxdt .
\]
where $\frac{\partial F(x,u)}{\partial u}=f(x,u)$, first seek weak solution of the modified equation:
\be
\Box u=\beta v_{tt}-f(x,u)-\beta v
\label{mpapier}
\ee
and then send the parameter $\beta$ to zero.\\
{\bf Assumptions on $f(u)$:}\\
we assume that there are positive constants $c_0^1 \leq c_0^2,c^1_1,c_1^2$ such that
\be
c_0^1|u|^{s-1}u+c_1^1\leq f(x,u) \leq c_0^2|u|^{s-1}u+c_2^1
\label{conditions}
\ee
with $c_0^1>\frac{c_0^2}{s+1}$. These assumptions are satisfied by some nonlinearities of polynomial type. $f(x,u)$ must also be strongly monotone increasing:
\be
\frac{\partial f(x,u)}{\partial u}\geq\alpha>0
\label{conditions2}
\ee

\begin{theo}Under assumptions (\ref{conditions}),(\ref{conditions2}) and $f\in C^{2,1}$,(\ref{papier}),(\ref{bdry}) admits infinitely many classical solutions.
\end{theo}
\section{Estimates}
Define $l^q=\{\hat{u}(j,k) s.t. \sum_{(j,k)\in\Z\times\Z}|\hat{u}(j,k)|^q <+ \infty\}$. 
\begin{theo}
\label{Riesz}The function $u=\sum_{2j\neq\pm k}\hat{u}(j,k)e^{2ijx+ikt}\in C^{\gamma}$ where $\gamma<1-\frac{1}{p}$ 
if
\be
\hat{u}(j,k)=\frac{\hat{f}(j,k)}{4j^2-k^2}
\ee
for $2j\neq \pm k$, 
$\hat{f}\in l^q$ and $\frac{1}{p}+\frac{1}{q}=1$.
\end{theo}
Proof:\\
Let $B_m$ the set
\[
B_m=\{ (j,k)\in\Z\times\Z \,\, 2|j|+|k|\leq 2.2^m \}
\]
and $\Delta_m$
\[
\Delta_m=B_m\setminus B_{m-1}
\]
so we have in $\Delta_m$
\[
2^m\leq 2|j|+|k|\leq 2.2^m
\]
and the $C^\gamma$ norm will be estimated by
\[
\sup_m 2^{\gamma m}||\Delta_m||_{C^0}\]
see \cite{Schlag} or \cite{Katznelson}. 
\begin{eqnarray}2^{\gamma m}||{\Delta}_m||_{C^0}
& = & 2^{\gamma m}||\sum_{(j,k)\in{\Delta}_m}\hat{u}(j,k)e^{i2jx}e^{ikt}||_{C^0}\nonumber\\
& = & ||\sum_{(j,k)\in\Delta_m}\frac{2^m\hat{f}(j,k)}{4j^2-k^2}e^{i2jx}e^{ikt}||_{C^0}\nonumber\\
& \leq & [\sum_{(j,k)\in\Delta_m}\frac{2^{\gamma mp}}{(|2j|^2-k^2)^p}]^{\frac{1}{p}}[\sum_{(j,k)\in \Z\times\Z}|\hat{f}(j,k)|^q]^{\frac{1}{q}}\nonumber\\
& \leq & [\sum_{(j,k)\in\Delta_m}\frac{(2|j|+|k|)^{\gamma p}}{(2|j|+|k|)^p(|2j|-|k|)^2)^p}]^{\frac{1}{p}}[\sum_{(j,k)\in \Z\times\Z}|\hat{f}(j,k)|^q]^{\frac{1}{q}}\nonumber\\
& \leq & [\sum_{(j,k)\in\Delta_m}\frac{1}{(2|j|+|k|)^{(1-\gamma) p}(|2j|-|k|))^p}]^{\frac{1}{p}}[\sum_{(j,k)\in \Z\times\Z}|\hat{f}(j,k)|^q]^{\frac{1}{q}}\nonumber\\
& \leq & c||\hat{f}||_{l^q}\leq c||f||_{L^p}
\end{eqnarray}
as long as $\gamma<1-\frac{1}{p}$ and the last inequality follows from the Hausdorff-Young theorem. \\
Remark: The argument here provides an alternate proof of the H\"older continuity of weak solutions of $\Box w=f$ where $f\in L^p\cap N^\perp$ where $N^\perp$ denotes the weak orthogonal of the kernel of $\Box$ with Dirichlet boundary conditions, proved by Brezis and Coron and Nirenberg in \cite{BCN} via Lovicarova's fundamental solution, for $1<p\leq 2$.\\
In the case that $p=2$ we have $u\in C^{0,\gamma}$ or similarly $f\in H^\alpha$ implies $u\in C^{\alpha+\frac{1}{2}}$.\\ Define
\be
u_{h_1,h_2}(x,t)=u(x+h_1,t+h_2)
\ee
and 
\be
\Delta_m^{++}=\{(j,k)\in \Z\times\Z \,(j,k)\in\Delta_m \,\, j\geq 0, k\geq 0\}
\ee
\be
\Delta_m^{+-}=\{(j,k)\in \Z\times\Z \, (j,k)\in\Delta_m \,\, j\geq 0, k<0\}
\ee
\be
\Delta_m^{-+}=\{ (j,k)\in \Z\times\Z \, (j,k)\in\Delta_m \,\, j<0, k\geq 0 \}
\ee
\be
\Delta_m^{--}=\{ (j,k)\in \Z\times\Z \, (j,k)\in\Delta_m \,\, j\leq 0, k< 0\}
\ee
and define $u^{--},u^{++},u^{-+},u^{+-}$ as 
\begin{eqnarray}\hat{u}^{++}(j,k) & = & \hat{u}(j,k) \,\ {\rm{if}} \,\ j\geq 0,k\geq 0\nonumber\\
                            & = & 0 \,\ {\rm{otherwise}}
\end{eqnarray}
\begin{eqnarray}\hat{u}^{+-}(j,k) & = & \hat{u}(j,k) \,\ {\rm{if}} \,\ j\geq 0,k< 0\nonumber\\
                            & = & 0 \,\ {\rm{otherwise}}
\end{eqnarray}
\begin{eqnarray}\hat{u}^{-+}(j,k) & = & \hat{u}(j,k) \,\ {\rm{if}} \,\ j< 0,k\geq 0\nonumber\\
                            & = & 0 \,\ {\rm{otherwise}}
\end{eqnarray}
\begin{eqnarray}\hat{u}^{--}(j,k) & = & \hat{u}(j,k) \,\ {\rm{if}} \,\ j< 0,k< 0\nonumber\\
                            & = & 0 \,\ {\rm{otherwise}}
\end{eqnarray}
\begin{lem} If $u^{++}\in C^{0,\gamma}$ then $u^{++}\in H^{\gamma^\prime}$ if $\gamma^\prime<\gamma$
\end{lem}
The analogue is also true for $u^{+-},u^{-+},u^{--}$.\\
Proof:\\
\begin{eqnarray}||u^{++}||_{H^{\gamma^\prime}}^2
& =    & \sum_m\sum_{(j,k)\in\Delta_m^{++}}(|2j|+|k|)^{2\gamma^\prime}|\hat{u}(j,k)|^2\nonumber\\
& \leq & \sum_m\sum_{(j,k)\in\Delta_m^{++}}2^{2(m+1)\gamma^\prime}|\hat{u}(j,k)|^2\nonumber\\
& \leq & \sum_m 2^{(m+1)\gamma^\prime}\sum_{(j,k)\in\Delta_m^{++}}|e^{i(2|j|+|k|)h}-1|^2|\hat{u}(j,k)|^2\nonumber
\end{eqnarray}
with $h=h(m)=\frac{2\pi}{3}2^{-m}$. \\
Remark: in the next line the sum over $\Delta_m^{++}$ is extended to the whole series but $h$ still depends on $m$. This is possible because of Parseval. \\Then
\begin{eqnarray}||u^{++}||_{H^{\gamma^\prime}}^2
& \leq & \sum_m2^{(m+1)\gamma^\prime}||u_{h(m),h(m)}-u||_{l^2}\\
\label{Uh}
& \leq & \sum_m2^{(m+1)\gamma^\prime}||u_{h(m),h(m)}-u||_{C^0}\nonumber\\
& \leq & \sum_m2^{(m+1)\gamma^\prime}||u||_{C^0,\gamma}^2|h(m)|^{2\gamma}\nonumber\\
& \leq & \sum_m2^{(m+1)\gamma^\prime}||u||_{C^0,\gamma}^2|\frac{2\pi}{3}2^{-m}|^{2\gamma}\nonumber\\
& \leq & c||u||_{C^0,\gamma}^2\sum_m2^{m(\gamma^\prime-\gamma)}\nonumber\\
& \leq & c(\gamma-\gamma^\prime)||u||^2_{C^{0,\gamma}}\nonumber
\end{eqnarray}
The estimates for $u^{--},u^{-+},u^{+-}$ follow similarly by replacing $u_{h,h}$ in the preceding argument (\ref{Uh}) by $u_{-h,-h}$,$u_{-h,h}$,$u_{h,-h}$.
We can conclude by noting:
\be
||u||^2_{H^\gamma}=||u^{++}||_{H^\gamma}^2+||u^{+-}||_{H^\gamma}^2+||u^{-+}||_{H^\gamma}^2+||u^{--}||_{H^\gamma}^2
\ee
We prove a bootstrapping estimate in the next lemma. It follows from the proof theorem 4 in \cite{R67} established for Dirichlet boundary conditions.
\begin{lem}
\label{H1}
Let $f,w\in L^2(Q)$ such that 
\be
\widehat{f}(j,k)=0=\widehat{w}(j,k)\,\ {\rm{for}}  \,\, 2j=\pm k
\ee
 and 
\be
(-k^2+4j^2)\widehat{w}(j,k)=\widehat{f}(j,k),\,\ {\rm{for}} \,\, 2j\neq\pm k \,\, 
\ee
then $w\in H^1$
\end{lem}
Proof:
\begin{eqnarray}||w||_{H^1}^2
& = & \sum_{2j\neq \pm k}\frac{4j^2+k^2}{|4j^2-k^2|^2}|\hat{f}(j,k)|^2\nonumber\\
& = & \sum_{2j\neq \pm k}\frac{1}{2}\frac{(2j-k)^2+(2j+k)^2}{(2j-k)^2(2j+k)^2}|\hat{f}(j,k)|^2\nonumber\\
& \leq & \sum_{2j\neq \pm k}|\hat{f}(j,k)|^2\nonumber\\
& \leq & ||f||_{L^2}^2.
\end{eqnarray}
Let $E^s$ be the closure of $\{e^{i2jx+ikt},\,\ 2j\neq \pm k\}$ under the norm
\[
||u||_{E^s}^2=\sum_{2j\neq\pm k}|\widehat{u}(j,k)|^2|k^2-4j^2|^s
\]
the we have the Sobolev estimate:
\begin{theo} $0<s<1$ the space $E^{s}$ is continuously embedded in $L^p$ where $p=\frac{2-s}{1-s}$.
\label{Sobolev}
\end{theo}
This theorem implies that the embedding in \cite{Coron} $E^1\subset L^p$ is compact, as $E^1\subset E^s$ is compact for $s<1$. We will show that it also implies a Gagliardo-Nirenberg inequality of the type:
\be
||u||_{L^p}\leq c(p)||u||_{L^2}^{1-s(p)}||u||_{E^1}^{s(p)}
\ee
where $c(p)$ will be computed explicitely.\\
Proof:\\
\be
f=f_{1,A}+f_{2,A}
\ee
where
\be
f_{1,A}=\sum_{2j\neq \pm k ,2|j|+|k|\leq A}\hat{f}(j,k)e^{i2jx}e^{ikt}
\ee
and
\be
f_{2,A}=\sum_{2j\neq \pm k ,2|j|+|k|>A}\hat{f}(j,k)e^{i2jx}e^{ikt}
\ee
\begin{eqnarray} |f_{1,A}| & \leq & \sum_{ 2j\neq \pm k ,2|j|+|k|\leq A}|\hat{f}(j,k)| \nonumber\\
                           & \leq & \sum_{2j\neq \pm k ,2|j|+|k|\leq A}|4j^2-k^2|^{-\frac{s}{2}}|4j^2-k^2|^{\frac{s}{2}}|\hat{f}(j,k)|
\end{eqnarray}
and applying Cauchy-Schwarz we have
\begin{eqnarray}|f_{1,A}|  & \leq & (\sum_{2j\neq \pm k ,2|j|+|k|\leq A}\frac{1}{|4j^2-k^2|^s})^{\frac{1}{2}}(\sum_{2j\neq \pm k ,2|j|+|k|\leq A}|4j^2-                                    k^2|^s|\hat{f}(j,k)|^2)^{\frac{1}{2}}\nonumber\\
                           & \leq & ||f||_{E^s}(\sum_{m,n\in\mathbb{N}\leq A}\frac{4}{m^sn^s})^{\frac{1}{2}}\nonumber\\
                           & \leq & c||f||_{E^s}(\int_1^A\frac{dm}{m^s}\int_1^{A}\frac{dn}{n^s})^{\frac{1}{2}}\nonumber\\
                           & \leq & c||f||_{E^s} A^{-s+1}.
\end{eqnarray}
Now we seek $A_\lambda$ such that 
\be
|f_{1,A}|\leq \frac{\lambda}{4}.
\label{chemin2}
\ee So we require the estimate
\be
c||f||_E A^{1-s}\leq \frac{\lambda}{4}
\label{chemin1}
\ee
this leads to the inequality
\[ A^{1-s}\leq \frac{\lambda}{4c||f||_{E^s}}.\]
So let $A_\lambda$:
\[ A_\lambda= (\frac{\lambda}{4c||f||_{E^s}})^{\frac{1}{1-s}}.  \]
Now 
\[ \int_{[0,\pi][0,2\pi]}|f(x,t)|^pdxdt=p\int_0^\infty y^{p-1}w(y)dy \]
where $w_f(y)=|\{(x,t)\in [0,\pi][0,2\pi]: |f(x,t|>y\}|$. Now $|f(x,t)|>\lambda$ implies $|f_{1,A}|>\frac{\lambda}{2}$ or
$|f_{2,A}|>\frac{\lambda}{2}$. Recalling (\ref{chemin2}) and the definition of $A_\lambda$ conclude that 
\be
|f_{2,A_\lambda}|>\frac{\lambda}{2}
\ee                       
and
\be
w_f(\lambda)\leq w_{f_{2,A_\lambda}}(\frac{\lambda}{2})
\ee
hence 
\begin{eqnarray}\int_{[0,\pi][0,2\pi]}|f(x,t)|^pdxdt & = & p\int_0^\infty \lambda^{p-1}w_f(\lambda)d\lambda \nonumber\\
                                                     & \leq & p\int_0^\infty \lambda^{p-1}w_{f_{2,A_\lambda}}(\frac{\lambda}{2})d\lambda\nonumber.
\end{eqnarray}
Since 
\be
w(\lambda)\leq \frac{1}{\lambda^2}\int_{|f|\geq \lambda}|f(x,t)|dxdt
\ee
\begin{eqnarray} \int_{[0,\pi][0,2\pi]}|f(x,t)|^pdxdt 
& \leq & \int_0^\infty\lambda^{p-3}\int_{|f_{2,A_\lambda}(x,t)|>\frac{\lambda}{2}}|f_{2,A_\lambda}(x,t)|^2dxdtd\lambda\nonumber\\
& \leq & \int_0^\infty\lambda^{p-3}\int_{[0,\pi][0,2\pi]}|f_{2,A_\lambda}(x,t)|^2dxdtd\lambda.
\end{eqnarray}
Then we can invoke Parseval formula to deduce
\begin{eqnarray} \int_{[0,\pi][0,2\pi]}|f(x,t)|^pdxdt 
& \leq & \int_0^\infty\lambda^{p-3}\sum_{2j\neq\pm k}|\hat{f}_{2,A_\lambda}(j,k)|^2d\lambda\nonumber\\
& = & \int_0^\infty\lambda^{p-3}\sum_{2j\neq \pm k,2|j|+|k|>A_{\lambda}}|\hat{f}(j,k)|^2d\lambda.
\label{f}
\end{eqnarray}
Now
\[ 2|j|+|k|\geq A_\lambda= (\frac{\lambda}{4c||f||_{E^s}})^{\frac{1}{1-s}} \] 
implies
\[ \lambda \leq 4c||f||_{E^s}((2|j|+|k|))^{1-s}\]
We continue the estimate from (\ref{f}):
\begin{eqnarray} \int_{[0,\pi][0,2\pi]}|f(x,t)|^pdxdtd\lambda 
& \leq & \sum_{2j\neq\pm k}\int_0^\infty\lambda^{p-3}|\hat{f}(j,k)|^21_{\{(\lambda,j,k) s.t. 2|j|+|k|\geq A_\lambda\}}d\lambda\nonumber\\
& \leq & \sum_{2j\neq\pm k}\int_0^{4c||f||_{E^s}(2|j|+|k|)^{1-s}}\lambda^{p-                                    3}|\hat{f}(j,k)|^2d\lambda\nonumber\\
& \leq & \sum_{2j\neq\pm k}|\hat{f}(j,k)|^2\int_0^{4c||f||_{E^s}(2|j|+|k|)^{1-s}}\lambda^{p-                                    3}d\lambda\nonumber\\
& \leq & \sum_{2j\neq\pm k,2|j|+|k|}|\hat{f}(j,k)|^2[\frac{\lambda^{p-2}}{p-2}]_0^{4c||f||_{E^s}(2|j|+|k|)^{1-s}}\nonumber\\
& \leq & \sum_{2j\neq\pm k,2|j|+|k|}|\hat{f}(j,k)|^2\frac{1}{p-2}[4c||f||_{E^s}((2|j|+|k|))^{1-s}]^{p-2}\nonumber\\
\end{eqnarray}
now if $s=(1-s)(p-2)$ i.e. $s(p)=\frac{p-2}{p-1}$ then
\be
\int_{[0,\pi][0,2\pi]}|f(x,t)|^pdxdt\leq \frac{(4c)^{p-2}}{p-2}||f||_{E^{s(p)}}^p
\ee
and we have the following Gagliardo-Nirenberg inequality for $p>2$:
\be
||u||_{L^p}\leq c(p)||u||_{E^{s(p)}}\leq c(p)||u||_{L^2}^{1-s(p)}||u||_{E^1}^{s(p)}.
\ee
\section{Construction of the weak solution}
 For the Galerkin procedure we define the spaces:
\[
E^m=span\{\sin 2jx\cos kt, \sin 2jx\sin kt,\, \cos 2jx\cos kt,\, \cos 2jx\sin kt,  \,\,\, 2j+k\leq m, 2j\neq k \},
\]
\[
E^{-m}=span\{\sin 2jx\cos kt, \sin 2jx\sin kt,\, \cos 2jx\cos kt,\, \cos 2jx\sin kt, \,\,\, 2j+k\leq m \,\ 2j<k \},
\]
\[
E^{+l}=span\{\sin 2jx\cos kt, \sin 2jx\sin kt,\, \cos 2jx\cos kt,\, \cos 2jx\sin kt, \,\,\, 2j+k\leq l \,\ 2j>k \},
\]
\[
N^m=span\{ \sin 2jx\cos kt, \sin 2jx\sin kt,\, \cos 2jx\cos kt,\, \cos 2jx\sin kt, \,\,\ 2j\leq m \}
\]
which are employed in the minimax procedure. We denote by $P^m$ the projection of $E\oplus N$ into $E^m\oplus N^m$. The functional $I_\beta$ satisfies the Palais-Smale condition. The arguments follows as in \cite{R84}, we do not repeat them here.
\begin{lem}$\forall u\in E^{+l}$, there is a constant $C(l)$ independent of $\beta,m$ such that 
\be
I_\beta(u)\leq M(l)
\ee
\end{lem} 
Proof:\\
Let $u\in E^{+l}$
\begin{eqnarray}I_\beta(u)
& =    & \frac{1}{2}||w^+||_E^2-\frac{1}{2}||w^-||_E^2-\beta||v||_{H^1}^2-\int_QF(u)dxdt\nonumber\\
& \leq & \frac{1}{2}||w^+||_E^2-\frac{1}{2}||w^-||_E^2-\beta||v||_{H^1}^2-c(s)\int_Q\frac{|u|^{s+1}}{s+1}dxdt+d(f,s)\nonumber\\
& \leq & c(f,s)+\sup_{u\in E^{+n}}\frac{1}{2}||w^+||_E^2-c(s,Q)||u||_{L^2}^{s+1}.
\end{eqnarray}
Now in $E^{+l}$
\be
||u||_E^2\leq l||u||_{L^2}^2
\ee
and on the other-hand
\be
\sup_{u\in E^{+l}}\frac{1}{2}||w^+||_E^2-c(s,Q)||u||_{L^2}^{s+1}>0
\ee
while as $||u||_E\rightarrow +\infty$ in $E^{+l}$ is dominated by $||u||_{L^2}^{s+1}$ as $s+1>2$ 
and is attained at say $\overline{u}$
hence we have 
\be
c(s,Q)||\overline{u}||_{L^2}^{s+1}\leq ||\overline{u}||_E^2\leq l||\overline{u}||_{L^2}^2
\ee
and we can conclude there is $M(l)$ depending on $l$ but independent of $\beta$ such that
\be
I_\beta(u)\leq M(l).
\ee
Also $E^{+l}$ is finite dimensional hence there is $R(l)$ such that for all $u\in E^{+l}\oplus E^{-m}\oplus N^m$ and $||u||_{E,\beta}\geq R(l)$ implies $I_{\beta}(u)\leq 0$.
\begin{theo} {\rm{Let}} $f$ be $C^1$, {\rm{for}} $l$ {\rm{large enough there is a distributional solution}} $u=v+w$ {\rm{of the modified problem}} (\ref{mpapier}).
\end{theo}
Proof:\\
In this proof the constants may dependent on $\beta$ and $f$ but are independent of $m$.
The proof of this theorem here is slightly simpler from the one in \cite{R84} as we take advantage of the polynomial growth of the nonlinear term. We also  employ Galerkin approximation.\\
Let $u_l^m=w^m+v^m\in E^m\oplus N^m$ a distributional solution corresponding to the critical value $c_l$, and any $\phi\in E^m\oplus N^m$:
\be
I^\prime(u_l^m)\phi=0
\label{Galerkin1}\ee
now taking $\phi=v_{tt}^m\in N^m$ we have
\[
(\beta v_{tt}^m,v_{tt}^m)_{L^2}=(f(x,u_l^m),v_{tt}^m)_{L^2}+\beta(v_t^m,v_t^m)
\]
and by (\ref{conditions}) there are constant positive $c,d$ such that 
\[
\beta||v_{tt}^m||_{L^2}^2\leq c||u^s||_{L^2}||v_{tt}^m||_{L^2}+d||v_{tt}^m||_{L^2}
\]
\[
\beta||v_{tt}^m||_{L^2}\leq c||v_{tt}^m||_{L^2}
\]
hence
\[
||v_{tt}^m||_{L^2}\leq c(\beta)
\]
we now have
\[
w_{tt}^m-w_{xx}^m=\beta v_{tt}^m+P^mf(x,u^m_l)\in L^2
\]
hence 
$w^m\in H^1\cap C^{\gamma}$, $\gamma<\frac{1}{2}$
by theorem \ref{Riesz} and lemma \ref{H1}. This now implies $w^m\in H^2$, $w^m\rightarrow w$ as $m\rightarrow +\infty$ pointwise and $w\in H^1\cap C^{\gamma}$. Then if $\phi=v_{tttt}^m$ then
\[
(\beta v_{tt}^m,v_{tttt}^m)_{L^2}=(f(x,u_l^m),v_{tttt}^m)_{L^2}-\beta(v^m_t,v^m_{ttt})
\]
so there exists $c$ independent of $m$ 
\[
(\beta v_{ttt}^m,v_{ttt}^m)_{L^2}= (f_u(x,u_l^m)u_{lt}^m,v_{ttt}^m)_{L^2}-\beta(v^m_t,v_{ttt}^m)
\] 
and we deduce $||v_{ttt}^m||_{L^2}\leq c(\beta)$ hence $v^m_{ttt}\rightarrow v_{tt}\in C^0$ as $m\rightarrow+\infty$ hence $v$ is $C^2$ and $w$ is $C^\gamma$ by applying theorem \ref{Riesz} to (\ref{mpapier}) . We now have
\[
u_l^m\rightarrow u\in C^\gamma \,\, {\rm{as}}\,\ m\rightarrow +\infty
\]
and since (\ref{Galerkin1}) holds for any $\phi\in E^m\oplus N^m$ we can deduce
\be
I^\prime(u)\phi=0 \,\,\, \forall \phi\in E^m\oplus N^m,
\ee
now sending $m\rightarrow\infty$,  $u$ is a weak solution of (\ref{mpapier}).\\
Then we can define $g_\theta(u)=u(x,t+\theta)$. Define:
\be
{\cal G}=\{ g_\theta \,\ s.t. \,\, \theta\in [0,2\pi) \}
\ee
\be
V_l=N^m\oplus E^{-m}\oplus E^{+l}
\ee
\be
G_l=\{h\in C(V_l,E^m) \,\, {\rm{ such \ that}} \,\ h \,\, {\rm{satisfies}} \,\,\ \gamma_1-\gamma_4 \}
\ee
$Fix {\cal G}=\{u\in E \,\, s.t. \,\, g(u)=u \,\, \forall g \in {\cal G} \}=span \{ \cos{2jx},\,\, \sin 2jx,\,\  j\in\Z \}\subset E^-$. Define $P^{0m},P^{-m}$ the orthogonal projections from $E^m\oplus N^m$ onto respectively $N^m(=E^{0m}),E^{-m}$ and $P_l$ the orthogonal projection from $E^m\oplus N^m$ onto $V_l$.
\begin{displaymath}
\left\{ \begin{array}{l}
\gamma_1 \,\, h \,\ {\rm{is \ equivariant}} \\
\gamma_2 \,\, h(u)=u \,\ {\rm{if}} \,\ u\in Fix{\cal G}      \\
\gamma_3 \,\ {\rm{There exists }} r=r(h) \,\ h(u)=u \,\, {\rm{if}} \,\, u\in V_l\setminus B_{r(h)} \\
\gamma_4 \,\ u=w^++w^-+v\in V_l \,\ (P^{0m}+P^{-m})h(u)=\alpha(u)v+\alpha^-(u)w^-+\phi(u) \,\, {\rm{where}} \,\, \alpha,\ov{\alpha}\in C(V_l,[1,\ov{\alpha}]) \\
\end{array}\right.
\end{displaymath}
and $1<\ov{\alpha}$ depends on $h,\phi$ continuous.
Define
\be
c_l(\beta)=\inf_{h\in G_l}\sup_{u\in V_l}I_\beta(h(u))
\ee
and $c_l(\beta)\rightarrow +\infty$ as $l\rightarrow\infty$ independently of $m,\beta$ .
\begin{lem}
$c_l(\beta)\rightarrow +\infty$ as $l\rightarrow +\infty$ 
\end{lem}
Proof:\\
\begin{eqnarray}I_\beta(u)
& = & \frac{1}{2}||w^+||_E^2-\frac{1}{2}||w^-||_E^2-\beta||v||_{H^1}^2-\int_Q F(u)dxdt
\end{eqnarray}
and there exists by assumptions (\ref{conditions}) $c(s),d(s)>0$ such that
\begin{eqnarray}I_\beta(u)
& \geq & \frac{1}{2}||w^+||_E^2-\frac{1}{2}||w^-||_E^2-\beta||v||_{H^1}^2-               c(s)\int_Q|u|^{s+1}dxdt-d(s)\nonumber\\
\end{eqnarray}
and if $u\in \partial B_\rho\cap V_{l-1}^\perp$ we have
\begin{eqnarray}I_\beta(u)
& \geq & \frac{1}{2}||w^+||_E^2-c(s)\int_Q|u|^{s+1}dxdt-d(s).\nonumber\\
\end{eqnarray}
Now by the Sobolev embedding theorem \ref{Sobolev} there is $\theta(s)<1$ such that if $\widehat{u}(j,k)=0$ for $2j=\pm k$ we have
\be
||u||_{L^{s+1}}\leq ||u||_{E^{\theta(s)}}
\ee
hence
\begin{eqnarray}I_\beta(u)
& \geq & \frac{1}{2}||w^+||_E^2-c(s)(||u||_{L^2}^{1-\theta(s)}||u||_{E^1}^{\theta(s)})^{s+1}-                  d(s)\nonumber\\
& \geq & \frac{1}{2}\rho^2-\rho^{s+1}l^{(1-\theta(s))\frac{s+1}{s-1}}.
\end{eqnarray}
If we choose a constant $C(s)$ large and $\rho=\frac{1}{C(s)}l^{(1-\theta(s))\frac{s+1}{s-1}}$
\begin{eqnarray}I_\beta(u)
& \geq & \frac{1}{4}\rho^2-d(s).\nonumber\\
\end{eqnarray}
Applying the corollary $2.4$ in \cite{FHR} to $P_{l-1}h\in C(\partial B_{\rho_l},V_{l-1})$
we have
\be
h(V_l)\cap\partial B_{\rho_l}\cap V_{l-1}^\perp\neq \o
\ee
hence
\be
\sup_{V_l}I_\beta(h(u))\geq \inf_{u\in\partial B_{\rho_l}\cap V_{l-1}^\perp}I(u)\geq \frac{1}{4}\rho^2-d(s)\rightarrow +\infty
\ee
The $c_l(\beta)$ are critical values of $I_\beta$ on $E^m$. This is obtained by a standard argument see \cite{R84} propositions 2.33 and 2.37. 
\begin{lem}If $u$ is a critical point of $I_\beta$ in $E^m\oplus N^m$ then there are constants $c_1,c_2$ independent of $m,\beta$ such that 
\be
||f(u)||_{L^{\frac{s+1}{s}}}^{\frac{s+1}{s}}\leq c_1I(u)+c_2
\ee
\end{lem}
Proof:\\
If $u$ is a critical point of $I_\beta$ then $I_\beta^\prime(u)\phi=0\,\, \forall \phi\in E^m\oplus N^m$
hence
\begin{eqnarray}I_\beta(u)
& = &     I_\beta(u)-I_\beta^\prime(u)u\nonumber\\
& = &     \int_Q\frac{1}{2}uf(u)-F(u)dxdt\geq a_1(s)\int_Q|u|^{s+1}dxdt-a_2(s)
\label{L2}
\end{eqnarray}
such constant $a_1(s),a_2(s)$ exist because $f$ satisfies (\ref{conditions}). Then we have
\begin{eqnarray}I_\beta(u)
& \geq & c_1\int_Q|f(u)|^{\frac{s+1}{s}}dxdt-c_2(s).
\end{eqnarray}
Let $u^m=w^m+v^m$ the approximate solution on $E^m\oplus N^m$ then
\be
\widehat{\Box w^m}(j,k)=\widehat{f(u^m)}(j,k)
\ee
 $\forall 2j\neq k \in E^m$, hence by lemma \ref{Riesz} and the Hausdorff-Young we have
\be
||w^m||_{C^\gamma}\leq c
\ee
with $c$ independent of $m,\beta$. Hence we can conclude that $w=\lim_{m\rightarrow+\infty}w^m\in C^\gamma$ for any $\gamma<1-\frac{s}{s+1}$.

In the following lemma we follow closely the method of \cite{R78} to get an a priori estimate on $||v||_{C^0}$ independently of $\beta$.
\begin{lem} There is a constant $c$ independent of $\beta$ such that
\be
||v(\beta)||_{C^0}\leq c
\ee
\end{lem}
Proof:\\
First note that by (\ref{L2}):
\be
||v(\beta)||_{L^2}^2\leq||u||_{L^{s+1}}^{s+1}\leq c(l)+a_2(s)
\ee
so we already have a $L^2$ a priori estimate on $v$ independently of $\beta$. The point of this lemma is then to prove a $C^0$ estimate. We will discuss two cases:\\
Case 1:\\ $||v(\beta)||_{C^0} \leq 8||v(\beta)||_{L^2}$. Then we have a $C^0$ estimate on $v(\beta)$ independently of $\beta$.\\
Case 2:\\  $||v(\beta)||_{C^0}> 8|Q|||v(\beta)||_{L^2}$.\\

Let $\phi\in N$ then we have
\be
\int_Q[-\beta v_{tt}+\beta v+f(v+w)]\phi dxdt=0
\ee
or
\be
\int_Q\beta v\phi+\beta v_{t}\phi_t+[(f(v+w)-f(w))\phi]dxdt=-\int_Qf(w)\phi dxdt
\ee
and $q$ is the function defined as 
\begin{displaymath}
q(s)=
\left\{ \begin{array}{ll}
s+M            & s\geq M \\
0           & -M\leq s\leq M \\
s-M & s<M \\
\end{array}\right.
\end{displaymath}
and choose 
\be
\phi(x,t)=q(v^+(x,t))+q(v^-(x,t)).
\ee
\begin{eqnarray}\int_Qv^-q^+dxdt
& = & \frac{1}{\frac{|Q|}{2}}\sum_{j,k}\widehat{v^-}(j,k)\overline{\widehat{q^+}(j,k)}\nonumber\\
& = & \frac{1}{{\frac{|Q|}{2}}}\widehat{v^-}(0,0)\overline{\widehat{q^+}(0,0)}\nonumber\\
& \leq & ||v^-||_{L^2}||q^+||_{L^1}=||v^-||_{L^2}||q^+||_{L^1}
\label{vq}
\end{eqnarray}
and $\int_Qv^+q^-dxdt=\frac{1}{\frac{|Q|}{2}}\widehat{v^+}(0,0)\overline{\widehat{q^-}(0,0)}=0$ similarly.
\begin{eqnarray}\int_Qv_t\phi_tdxdt
& = & \int_Qq^\prime(v^+)(v^{+}_t)^2+q^\prime(v^-)(v^{-}_t)^2+\frac{\partial}{\partial_t}(q(v^+))v^-_t+\frac{\partial}{\partial_t}(q(v^-))v^+_tdxdt\nonumber\\
& = & \int_Qq^\prime(v^+)(v^{+}_t)^2+q^\prime(v^-)(v^{-}_t)^2dxdt,
\end{eqnarray}
we define
\begin{displaymath}
\psi(z)=
\left\{ \begin{array}{ll}
\min_{|\xi|\leq M_5}f(z+\xi)-f(\xi)       & z\geq 0 \\
\max_{|\xi|\leq M_5}f(z+\xi)-f(\xi)       & z<0     \\
\end{array}\right.
\end{displaymath}
which is monotone in $z$ with $\psi(0)=0$.
$Q_\delta=\{(x,t)\in Q,\,  |v(x,t)|\geq \delta \}$, $Q_{\delta}^+=\{(x,t)\in Q \,\, v(x,t)\geq \delta \}$, $Q_\delta^-=Q_\delta\setminus Q_\delta^+$.
If $v\geq 0$ then 
\begin{eqnarray}\int_{Q_\delta^+}[f(v+w)-f(w)][q^++q^-]dxdt
& \geq & \frac{\psi(\delta)}{||v|_{C^0}}\int_{{Q_\delta}^+}v(q^++q^-)dxdt.\nonumber
\end{eqnarray}
Now for $v\leq -\delta$
when $v<0$ then $f(v+w)-f(w)\leq \psi(v)$ and also $q^++q^-\leq 0$
and similarly
\begin{eqnarray}\int_{Q_\delta^-}(f(v+w)-f(w)(q^++q^-)dxdt
& \geq & \frac{-\psi(-\delta)}{||v||_{C^0}}\int_{{Q_\delta}^-}v(q^++q^-)dxdt\nonumber
\end{eqnarray}
now define $\nu(z)=\min (\psi(z),\psi(-z))$ for $z\geq 0$, and $||v^\pm||_{C^0}=\max(||v^+||_{C^0},||v^-||_{C^0})$
then
\begin{eqnarray}\int_{Q_\delta}[f(v+w)-f(w)][q^++q^-]dxdt
& \geq & \frac{\nu(\delta)}{||v||_{C^0}}[\int_{Q}v^+q^++v^-q^-dxdt-\delta\int_Q|q^+|+|q^-|dxdt]\nonumber\\
&      & -\frac{\nu(\delta)}{||v||_{C^0}}||v^-||_{L^2}\int_Q|q^+|+|q^-|dxdt\nonumber
\end{eqnarray}
where we employ (\ref{vq}) to estimate the terms $\int_Qv^+q^-+v^-q^+dxdt$,
and since $sq(s)\geq M|q(s)|$, then
\be
||f(w||_{C^0}\int_Q(|q^+|+|q^-|)dxdt\geq \frac{(M-||v^-||_{L^2}-\delta)\nu(\delta)}{||v||_{C^0}}\int_Q(|q^+|+|q^-|)dxdt
\ee
then for arbitrary $M<|||v^\pm||_{C^0}$ and choosing $\delta=\frac{|||v^\pm||_{C^0}}{2}$ we deduce
\be
\nu(\frac{1}{2}||v^\pm||_{C^0})\leq 8||f(w||_{C^0}
\ee
hence $||v||_{C^0}$ is bounded independently of $\beta$. \\
\section{Regularity of the solution}
Here we prove that if $f\in C^{2,1}$ then the weak solution $u$ is $C^2$. Since $||v||_{C^0},||w||_{C^0}$ are bounded independetly of $\beta$, $f(u)\in C^0$. We also have
\be
(-4j^2+k^2)\widehat{w}(j,k)=\widehat{f(x,v+w)}(j,k)\,\,\ 2j\neq\pm k.
\label{perp}
\ee
Then by lemma \ref{H1} we have $w\in H^1$. Since $f$ is smooth then too $f(w+v)\in H^1$. Then (\ref{perp}) implies $w\in H^2$ and iterating once again leads to $w\in H^3$. Now going back to the original equation:
\be
-\beta v_{tt}=\Box w-f(x,u)-\beta v
\ee
and recalling that $v\in C^2$ we deduce  $v\in H^3$ which with (\ref{perp}) implies $w\in H^4$. Iterating once more implies $v\in H^4$ then again $w\in H^5$ and $v\in H^5$.
 Thus we can differentiating with refer to $t$ in the weak sense and we have
\be
\Box w_t-\beta v_{ttt}=-f_u^\prime(x,u)(w_t+v_t)-\beta v_t
\label{perp3}
\ee
in Fourier space. We now want to get estimates independently of $\beta$ pass to the limit and find solutions of (\ref{papier}). 
Now multiplying by $v_t(\beta)$ (this is possible
 since $v\in H^1$) and integrating we have
\be
\beta(v_t,v_t)+\beta (v_{tt},v_{tt})+(f_u^\prime(x,u)v_t,v_t)=-(f_u^\prime(x,u)w_t,v_t)
\label{perp2}
\ee
and
\be
\alpha ||v_t||_{L^2}^2<(f_u^\prime(x,u)v_t,v_t) \leq -(f_u^\prime(x,u)w_t,v_t)
\ee
Now since $f_u^\prime>\alpha>0$ and $||w||_{H^1}\leq c$ with $c$ independently of $\beta$ hence there is a constant $c$ independent of $\beta$
such that $||v_t||_{L^2}\leq c$. This combined with (\ref{perp}) implies $||w||_{H^2}\leq c$ where $c$ is independent of $\beta$. Differentiating (\ref{perp3}) with refer to $t$ we get 
\be
\beta v_{tt}+\Box w_{tt}-\beta v_{tttt}+f_u^\prime(x,u)v_{tt}=-f_{uu}^{\prime\prime}(x,u)v_t^2-f_{uu}^{\prime\prime}(x,u)w_t^2-2f_{uu}^{\prime\prime}(x,u)w_tv_t-f_u^\prime(x,u)w_{tt}
\label{perp6}
\ee
Now we multiply (\ref{perp6}) by $v_{tt}$ and estimate the $L^2$ norm of the first term of the RHS. 
\begin{eqnarray} (f_{uu}^{\prime\prime}(x,u)v_t^2,v_{tt})
& \leq & c(f)\int_0^\pi\int_0^{2\pi}v_t^2|v_{tt}|dxdt\nonumber\\
& \leq & c(f)(\int_0^\pi\int_0^{2\pi}v_t^4dxdt)^{\frac{1}{2}}(\int_0^\pi\int_0^{2\pi}|v_{tt}|^2dxdt)\nonumber
\end{eqnarray} 
we then deduce
\begin{eqnarray} (f_{uu}^{\prime\prime}(x,u)v_t^2,v_{tt})
& \leq & c(f)||v_t||_{L^2}^{\frac{3}{2}}||v_t||_{H^1}^{\frac{1}{2}}\\
\label{GagNC}
& \leq & c(f)||v_t||_{H^1}^{\frac{3}{2}}\nonumber
\label{perp7}
\end{eqnarray}
where the constant $c(f)$ is independent of $\beta$ and the inequalities in the previous argument stems from the Gagliardo-Nirenberg inequality.

The $L^1$ norms of the terms in the RHS of (\ref{perp6}) multiplied by $v_{tt}$ can be estimated by noting that $f(u)\in H^1$,$w\in C^{1,\gamma}$,$0<\gamma<\frac{1}{2}$, and that the respective norms can be estimated are independently of $\beta$:
\be
(f_{uu}^{\prime\prime}(x,u)w_t^2,v_{tt}) \leq c||v_{tt}||_{L^2}
\label{perp8}
\ee
\be
(2f_{uu}^{\prime\prime}(x,u)w_t,v_{tt}) \leq c||v_{tt}||_{L^2}
\label{perp9}
\ee
\be
-(f_u^\prime(x,u)w_{tt},v_{tt})\leq c||w_{tt}||_{L^2}||v_{tt}||_{L^2}
\label{perp10}
\ee
recalling (\ref{perp6}), multiplying by $v_{tt}$
\be
\beta(v_{tt},v_{tt})+\beta(v_{ttt},v_{ttt})+(f_u^\prime(x,u)v_{tt},v_{tt})=(-f_{uu}^{\prime\prime}(x,u)v_t,v_{tt})+(-f_{uu}^{\prime\prime}(x,u)w_t,v_{tt})
+(-f_u^\prime(x,u)w_{tt},v_{tt})
\label{perp5}
\ee
We can now continue from (\ref{perp6}),(\ref{perp7}),(\ref{perp8}),(\ref{perp9}),(\ref{perp10}) and we have
\begin{eqnarray}\beta(v_{tt},v_{tt})+\beta(v_{ttt},v_{ttt})+(f_u^\prime(x,u)v_{tt},v_{tt})
& \leq & c||v_t||_{H^1}^{\frac{3}{2}}
\end{eqnarray}
thus there exists $c$ independent of $\beta$ such that $||v_{tt}||_{L^2}\leq c$ where $c$ is independent of $\beta$. At this stage we can conclude that there is a constant $c$ independent of $\beta$ such that $||f(u)||_{H^2}\leq c$. Combining this with (\ref{perp}) we have $||w||_{H^3}\leq c$ with $c$ is independent of $\beta$, $w\in C^{2,\frac{1}{2}}$ and $v\in C^1$ with upper bounds independent of $\beta$. We have now proved that if $f$ is $C^2$ then the solution is $u\in H^2\cap C^1$ is a weak solution of the equation.
We now differentiate (\ref{perp6}) we have
\begin{eqnarray}\beta v_{ttt}+\Box w_{ttt}+f^\prime(u)v_{ttt}
& = & -f_{uu}^{\prime\prime}(x,u)v_{tt}-f_{uuu}^{\prime\prime\prime}(x,u)(v_t+w_t)v_t^2-f_{uu}^{\prime\prime}(x,u)2v_tw_{tt}-f_{uuu}^{\prime\prime\prime}(x,u)w_{t}^2\nonumber\\
&   & -f_{uu}^{\prime\prime}(x,u)2w_tw_{tt}-2f_{uuu}^{\prime\prime\prime}(x,u)(v_t+w_{tt})w_{t}v_t-2f_{uu}^{\prime\prime}(x,u)w_{tt}v_t\nonumber\\
&   & -2f_{uu}^{\prime\prime}(x,u)w_tv_{tt}-f_{uu}^{\prime\prime}(x,u)(v_t+w_t)w_{tt}-f_{uu}^{\prime\prime}(x,u)w_{ttt}\nonumber
\end{eqnarray}
and multiplying both sides of the preceding equality by $v_{ttt}$ and integrating we conclude that $||v_{ttt}||_{L^2}\leq c$ where $c$ is independent of $\beta$ thus $v$ is $C^2$. Now recalling the Holder regularity bootstrap and (\ref{perp}) we get $w\in C^{3,\gamma}$, $0<\gamma<\frac{1}{2}$.
\bibliographystyle{plain}

\end{document}